\def\draft{n}
\documentclass{amsart}
\usepackage{
  fullpage,amssymb,epic,eepic,epsfig,verbatim}

\theoremstyle{plain}

\newtheorem{theorem}{Theorem}
\newtheorem{proposition}{Proposition}[section]
\newtheorem{lemma}[proposition]{Lemma}
\newtheorem{corollary}[proposition]{Corollary}

\theoremstyle{definition}
\newtheorem{definition}[proposition]{Definition}

\theoremstyle{remark}
\newtheorem{example}[proposition]{Example}

\newtheorem{remark}[proposition]{Remark}

\def\printname#1{
	\if\draft y
		\smash{\makebox[0pt]{\hspace{-0.5in}
			\raisebox{8pt}{\tt\tiny #1}}}
	\fi
}

\newcommand{\psdraw}[2]
         {\begin{array}{c} \hspace{-1.3mm}
	\raisebox{-4pt}{\epsfig{figure=draws/#1.eps,width=#2}}
	\hspace{-1.9mm}\end{array}}

\newlength{\standardunitlength}
\catcode`\@=11
\long\def\@makecaption#1#2{%
     \vskip 10pt

\setbox\@tempboxa\hbox{
       \small\sf{\bfcaptionfont #1. }\ignorespaces #2}%
     \ifdim \wd\@tempboxa >\captionwidth {%
         \rightskip=\@captionmargin\leftskip=\@captionmargin
         \unhbox\@tempboxa\par}%
       \else
         \hbox to\hsize{\hfil\box\@tempboxa\hfil}%
     \fi}
\font\bfcaptionfont=cmssbx10 scaled \magstephalf
\newdimen\@captionmargin\@captionmargin=2\parindent
\newdimen\captionwidth\captionwidth=\hsize
\catcode`\@=12

\def\lbl#1{\label{#1}\printname{#1}}

\def\BQ{\mathbb Q}

\def\A{\mathcal A}

\def\K{\mathcal K}
\def\cL{\mathcal L}

\def\E{\mathcal E}

\def\l{\lambda}

\def\s{\sigma}

\def\w{\omega}

\def\e{\epsilon}

\def\d{\delta}

\def\s{\sigma}

\def\matt#1#2#3#4#5#6#7#8#9{\left[
\begin{matrix}
 #1 & #2 & #3 \\
 #4 & #5 & #6 \\
 #7 & #8 & #9   
\end{matrix}
\right]}

\def\Per{\mathrm{Per}}
\def\fsl{\mathfrak s\mathfrak l}
\def\LIG{\mathrm{LIG}}
\def\LID{\mathrm{LID}}
\def\TID{\mathrm{TID}}
\def\IM{\mathrm{IM}}
\def\IMJ{\mathrm{IM}_J}
\def\ot{\otimes}

\begin{document}


\title[Random walks and the colored Jones function]{Random walks and the 
colored Jones function}

\author{Stavros Garoufalidis}
\address{Department of Mathematics \\
          University of Warwick \\
          Coventry, CV4 7AL, UK. }
\email{stavros@maths.warwick.ac.uk}
\author{Martin Loebl}
\address{KAM MFF UK and ITI \\
Charles University \\
Malostranske n. 25 \\
118 00 Praha 1 \\
Czech Republic.}
\email{loebl@kam.ms.mff.cuni.cz}

\thanks{S.G. was partially supported by an NSF and by an Israel-US BSF grant. 
M.L. was partly supported by GAUK 158 grant and by the Project LN00A056 of 
the Czech Ministry of Education.
        This and related preprints can also be obtained at
{\tt http://www.math.gatech.edu/$\sim$stavros } 
\newline
1991 {\em Mathematics Classification.} Primary 57N10. Secondary 57M25.
\newline
{\em Key words and phrases: Colored Jones function, permanents, weight 
systems, random walks.} 
}

\date{
This edition: February 28, 2002 \hspace{0.5cm} First edition: February 28, 
2002.}


\begin{abstract}
It can be conjectured that the colored Jones function of a knot can be
computed in terms of counting paths on the graph of a planar projection
of a knot. On the combinatorial level, the colored Jones function can
be replaced by its weight system. We give two curious formulas for the weight 
system of a colored Jones function: one in terms of the permanent of a matrix 
associated to a chord diagram, and another in terms of counting paths of 
intersecting chords.
\end{abstract}

\maketitle



\section{Introduction}
\lbl{sec.intro}
\subsection{History} 
\lbl{sub.history}

The colored Jones function of a knot is a $2$-parameter formal power series 
$\sum_{n,m=0}^\infty a_{n,m} h^n \l^m$ which determines (and is determined by) 
the Jones polynomial of a knot and its cables, \cite{BG}. The support of the
colored Jones function lies in the triangle $0 \leq m \leq n$. 

About 10 years ago, Melvin-Morton and Rozansky independently conjectured
a relation among the diagonal terms $\sum_n a_{n,n}(h\l)^n$ of the colored 
Jones function 
of a knot and its Alexander polynomial, \cite{MM,Ro1,Ro2}. D. Bar-Natan and 
the first author reduced the conjecture about knot invariants to a statement 
about their combinatorial weight systems, and then proved it for all 
semisimple Lie algebras using combinatorial methods, \cite{BG}.

Over the years, the MMR Conjecture has received attention by many researchers 
who gave alternative proofs, \cite{Ch,KSA,KM,Ro3,V}.

The subdiagonal terms $h^k\sum_n a_{n+k,n} (h \l)^n$ for a fixed $k$ of
the colored Jones function, are (after a suitable parametrization) rational
functions whose denominators are powers of the Alexander polynomial. This
was first shown by Rozansky in \cite{Ro3}, who further conjectured that a 
similar property should be hold for the full Kontsevich integral of a knot. 
Rozansky's conjecture was recently settled by the first author and Kricker 
in \cite{GK}. This opens the possibility of understanding each subdiagonal 
term of the colored Jones function (or the full Kontsevich integral) in 
topological terms.
That said, not much is known about the subdiagonal terms of the colored
Jones function. One can conjecture that each subdiagonal term
is given in terms of a certain counting of random walks on a planar 
projection of a knot, see also Lin and Wang, \cite{LW}.
On a combinatorial level, the colored Jones function may be replaced by
its weight system. In \cite{BG} formulas for the 
weight system $W_J$ of the colored Jones function and of its leading order 
term $W_{JJ}$ were given in terms of the intersection matrix of a chord 
diagram. In particular, $W_{JJ}$ equals to the {\em permanent} of the 
intersection matrix of the chord diagram. In the last section of \cite{BG} it 
was asked for a better understanding of the $W_J$ weight system, especially 
one that offers control over the subdiagonal terms in $W_J$.

The purpose of the paper is to give two curious combinatorial formulas for 
$W_J$ in Theorems \ref{thm.1} and \ref{thm.2} that answer these questions,
and support the conjecture that the colored Jones polynomial is a counting of
random walks.

\subsection{Statement of the results}
\lbl{sub.statement}

Consider the $0$-{\em framed colored Jones weight system}
$$
W_J :\A\to\BQ[\l]
$$
where $\A$ is the vector space over $\BQ$ spanned by {\em chord diagrams}
on an oriented line, modulo the $4$-term and $1$-term relations, see \cite{BG}
and also below. We will normalize $W_J$ to equal $1$ on the chord diagram
with no chords (in \cite{BG} the value of the empty chord diagram was
$\l+1$ instead). With this normalization, it turns out that for a chord 
diagram $D$,
$W_J(D)$ is a polynomial of $\l$ of degree the number of chords of $D$. 
$W_{JJ}(D)$ is defined to be the coefficient of $\l^{\mathrm{deg}}$ in 
$W_J$.
 
Given a chord diagram $D$, its chords are ordered (from left to right) and 
we can consider its intersection matrix $\IM(D)$ as in 
\cite[Definition 3.4]{BG} of size the number of chords of $D$ defined by
$$
\IM(D)_{ij}=\begin{cases}
\mathrm{sign}(i-j) & \text{if the chords $i$ and $j$ of $D$ intersect} \\
0 & \text{otherwise}.
\end{cases}
$$
We will consider a blown-up variant $\IMJ$ of the intersection matrix, of size
$3$ times the number of chords of $D$ composed of blocks of $3$ by $3$ matrices
as follows:
$$
\IMJ(D)_{ij}=\begin{cases}
A_{\mathrm{sign}(i-j)} & \text{if the distinct
         chords $i$ and $j$ of $D$ intersect} \\
A_0                    & \text{if $i=j$} \\
A_c & \text{if chords $i,j$ do not intersect and $i$ is completely
contained in $j$} \\ 
0 & \text{otherwise},
\end{cases}
$$
where
$$
A_0= \matt {\l+2} 0 0 0 {\l+2} 1 {\l} {-\l-2} 1 ,
\hspace{0.7cm}
A_-=\matt 0 0 0 {-1} {-1} 0 0 0 0 ,
\hspace{0.7cm}
A_+=\matt 1 1 0 0 0 0 0 0 0 , 
\hspace{0.7cm}
A_c=\matt 1 1 0 {-1} {-1} 0 0 0 0.
$$

\begin{example}
\lbl{ex.1}
$$ 
D=\psdraw{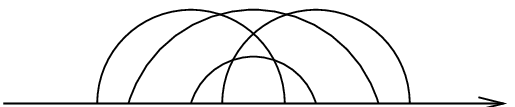}{1in}, \hspace{0.5cm} \IM(D)=
\left[
\begin{matrix}
0 & -1 & -1 & -1   \\
1 & 0 & 0   & -1 \\
1 & 0   & 0 & -1 \\
1   & 1 & 1 & 0
\end{matrix}\right], \hspace{0.5cm}
\IMJ(D)=
\left[
\begin{matrix}
A_0 & A_- & A_- & A_-   \\
A_+ & A_0 & 0   & A_- \\
A_+ & A_c   & A_0 & A_- \\
A_+ & A_+ & A_+ & A_0
\end{matrix}\right].
$$
\end{example}

\begin{theorem}
\lbl{thm.1}
We have
$$
W_J=\Per(\IMJ)
$$
where $\Per(A)$ denotes the permanent of a square matrix $A$.
\end{theorem}

There is an alternative (and equivalent) formula of $W_J$ in terms of 
counting cycles. In order to state it, given a chord diagram $D$ 
consider its {\em labeled intersection graph} 
$\LIG(D)$ as in \cite[Definition 3.4]{BG}. The vertices of $\LIG(D)$
correspond to the chords of $D$ (thus, are ordered) and the edges of $\LIG(D)$
correspond to the intersection of the chords of $D$.   

We will use a variation $\LID(D)$, the {\em labeled intersection digraph} of 
$D$ defined as follows. Orient each edge from the smaller vertex to the larger
and add an oriented loop on each vertex. The oriented loops are {\em leaving} 
the vertices. Next add directed edges $(ij)$ for each pair of chords
$i,j$ such that $i$ is completely contained in $j$. In addition we color these
new arcs {\em red} (and we draw them as \makebox[0.8cm]{$\longrightarrow
\hspace{-0.5cm}r$}) to 
distinguish them from the original arcs.

\begin{example}
\lbl{ex.22}
For the chord diagram $D$ of Example \ref{ex.1}, we have
$$
\LIG(D)=\psdraw{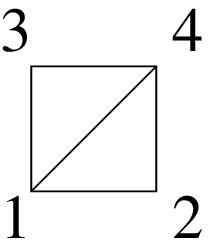}{0.6in}, \hspace{1cm}
\LID(D)=\psdraw{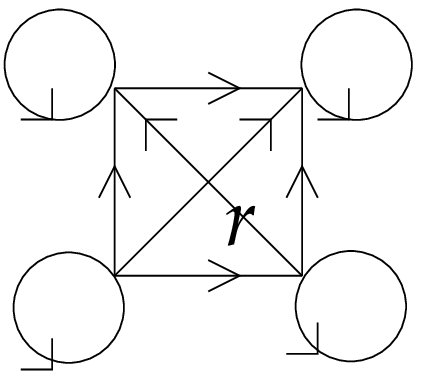}{0.7in}.
$$
\end{example}

A bit more generally, consider a {\em digraph} $G=(V,A)$ (i.e.,
a directed graph) where $V$ is the set
of vertices and $A$ is the set of arcs. If $e$ is an arc of $A$ with initial
vertex $u$ and terminal vertex $v$ then we write $e=(u,v)$.
 We assume that there is one loop at each vertex and
 a loop at a vertex is considered as an arc {\em leaving} that vertex,
and in addition some arcs which are not loops are {\em red}.
We will consider the arcs with variables associated with them: the variable
of an arc $e$ is denoted by $x_e$.  We will need the following notion
of {\em acceptable} object, given $G$:

\begin{definition}
\lbl{def.acceptable}
A collection $K$ of arcs together with a {\em thickening}
of one end of each of the arcs of $K$ is called {\em acceptable} for $G$
 if the following properties are satisfied:
\begin{itemize}
\item
If $(ij)$ is a red arc of $G$ then both arcs $(ij)$ and $(ji)$ may
appear in $K$, but they must always be thickened at $i$. If $(ij)$
is an uncolored arc of $G$ then  $(ij)$ with any end thickened
may appear in $K$, but $(ji)$ may not.
\item
Each vertex of $V$ is incident with $0$, $2$ or $4$ thickened arcs of $K$. If 
a loop belongs to $K$ then we assume it contributes 2 to the degree of the 
corresponding vertex. Moreover, a loop is always thickened at its initial 
segment, i.e., in agreement with its orientation.
\item
Exactly half of the arcs incident with a vertex are thickened at the vertex.
\item
If there are two arcs thickened at a vetex, then one of them enters and 
the other one leaves.
\end{itemize}
\end{definition}

We will study the following partition
function
$$
J(G)=\sum_{K \text{ acceptable}} 2^{\deg_4(K)} (\l+2)^{|V|-\deg_4(K)}
x_K (-1)^{a(K)}
$$
of a digraph $G$, 
where $x_K=\prod_{e \in K} x_e$, $\deg_4(K)$ denotes the number of vertices
of $K$ incident with 4 arcs of $K$ and $a(K)$ is the number 
of arcs of $K$ with initial segment thickened, i.e., directed in agreement 
with the thickening.

The motivation for $J(G)$ comes from the case of the intersection digraph
$\LID(D)$ of a chord diagram and the following:

\begin{theorem}
\lbl{thm.2}
For a chord diagram $D$, we have 
$$
W_J(D)=J(\LID(D))|_{x_e=1}.
$$
\end{theorem}


\begin{corollary}
\lbl{cor.2}
After a change of variables $d=\l+2$, let $W_{JJ^{(n)}}$ denote the coefficient
of $d^{\deg-n}$ in $W_J$.
Then,
$$
W_{JJ^{(n)}}(D)= 2^n\sum_{K} (-1)^{a(K)}
$$
where the sum is over all acceptable $K$ such that $\deg_4(K)=n$.
\end{corollary}

\begin{corollary}
\lbl{cor.3}\cite{BG} We have:
$$
W_{JJ}=\Per(\IM).
$$
\end{corollary}

How fast can one compute permanents? 
\begin{corollary}
\lbl{cor.4}
For general matrices of size $n$
we need $n!$ steps. However, a theorem of A. Galluccio \cite{GL} and the 
second author implies that $W_J$ can be computed in $4^g$ steps, where $g$ is 
the {\em genus} of $\LIG(D)$, that is the smallest genus of a surface that
$\LIG(D)$ embeds. 
\end{corollary}

\subsection{Plan of the proof}
\lbl{sub.plan}

In Section \ref{sec.review}, we review the weight system of the colored Jones
function, and reduce Theorem \ref{thm.2} to a Theorem \ref{thm.3} concerning
digraphs. Section \ref{sec.understanding}, is devoted to the proof of
Theorem \ref{thm.3} using a trip to combinatorics. 
In Section \ref{sec.per}, we translate our results using the language of
permanents, and deduce Theorem \ref{thm.1}.
In the final Section \ref{sec.cor} we prove the corollaries that follow
Theorem \ref{thm.2}.

\subsection{Acknowledgement}
The authors wish to thank the Georgia Institute of Technology which invited 
the second author and provided the environment of this research, and D. 
Bar-Natan for providing independent checks to the output of the program.

\tableofcontents

\section{A review of the $W_J$ weight system}
\lbl{sec.review}

The goal of this section is to reduce Theorem \ref{thm.2} to Theorem 
\ref{thm.3} stated below; this will be achieved 
by a careful examination of the $W_J$ weight system. Recall from 
\cite[Section 4.2]{BG} that $W_J$ can be computed as follows: \newline
{\bf Step~1.} Color each chord of a chord diagram $D$ by the following
operator:
\begin{eqnarray*}
\hat B(v_k \ot v_{k'}) & = &
(k+1)(\l-k'+1) v_{k+1} \ot v_{k'-1} \\
& + & (\l-k+1)(k'+1) v_{k-1} \ot v_{k'+1} \\
& + & 1/2((\l-2k)(\l-2k')-\l(\l+2)) v_k \ot v_{k'}
\end{eqnarray*}
from \cite[p.121]{BG}. \newline
The key calculation is the following elementary rearrangement of $\hat B$,
easily checked:

\begin{lemma}
\lbl{lem.1}
$$
\hat B(v_k \ot v_{k'})= ((\l+2) I + B^+ + B^-)(v_k \ot v_{k'})
$$
where
$$
B^+=\sum_{\e =0,1} (-1)^{\e} B^+_{\e}, \hspace{1cm}
B^+_{\e}=-(1+k)(\l+1-k') v_{k+\e} \ot v_{k'-\e} $$
$$
B^-=\sum_{\e =0,1} (-1)^{\e} B^-_{\e}, \hspace{1cm}
B^-_{\e}=-(1+k')(\l+1-k) v_{k-\e} \ot v_{k'+\e}.
$$
\end{lemma}

This done, the coloring of chords of $D$ may be viewed as a function
 $\rho:\text{chords}(D)\to \{I, B^+_0,B^+_1,B^-_0,B^-_1\}$. 

{\bf Step~2.} The end-points of the $n$ chords of $D$ partition
the base line into $2n+1$ segments $s_0,...,s_{2n}$ listed from left to right.
We associate number $m(s_i)$ with each of these segments as follows:
\begin{enumerate}
\item[1.] Let $m(s_0)=0$.
\item[2.] If $i\geq 0$ and last point of $s_i$ is left end-vertex of chord $v$
then $m(s_{i+1})$ is computed from $m(s_i)$ and $\rho(v)$ using \ref{lem.1}:
\newline
$\bullet$ If $\rho(v)\in \{I,B^+_0, B^-_0\}$ then $m(s_{i+1})=m(s_i)$, \newline
$\bullet$ If $\rho(v)=B^+_1$ then $m(s_{i+1})=m(s_i)+1$, \newline
$\bullet$ If $\rho(v)=B^-_1$ then $m(s_{i+1})=m(s_i)-1$. 
\item[3.] If $i\geq 0$ and last point of $s_i$ is right end-vertex of chord 
$v$ then $m(s_{i+1})$ is computed from $m(s_i)$ and $\rho(v)$ using 
\ref{lem.1}:\newline
$\bullet$ If $\rho(v)\in \{I,B^+_0, B^-_0\}$ then $m(s_{i+1})=m(s_i)$, \newline
$\bullet$ If $\rho(v)=B^+_1$ then $m(s_{i+1})=m(s_i)-1$, \newline
$\bullet$ If $\rho(v)=B^-_1$ then $m(s_{i+1})=m(s_i)+1$.
\end{enumerate}

{\bf Step~3.} We let
$$
W_J(D)=\sum_{\rho} \prod_{\text{chords}\, v} \w_\rho(v)
$$
where $\rho$ is a coloring of the chords of $D$ and $\w_\rho(v)$ is a specific
weight that is computed using Lemma \ref{lem.1} again:\newline
$\bullet$ $\w_\rho(v) = \l+2$ if $\rho(v) = I$,\newline
$\bullet$ $\w_\rho(v) = -(-1)^{\e}(1+k)(\l+1-k')$ if $\rho(v) = B^+_{\e}$,\newline
$\bullet$ $\w_\rho(v) = -(-1)^{\e}(1+k')(\l+1-k)$ if $\rho(v) = B^-_{\e}$,\newline
where $k=m(s_L(v))$, $k'= m(s_R(v))$, $s_L(v)$ is the segment ending at the 
left end-point of $v$ and $s_R(v)$ is the segment ending at the right 
end-point of $v$.
\medskip

The reader is urged to look at \cite[Chapter 4]{BG} for an explanation
of the above algorithm in terms of the representation theory of the 
$\fsl_2$ Lie algebra.

\begin{example}
\lbl{ex.3}
For the following coloring of the chord diagram of Example \ref{ex.1}
$$
\printname{bigCD}
	\setlength{\unitlength}{0.04\standardunitlength}
	\begin{array}{c}  \hspace{-1.7mm}
         	\raisebox{-8pt}{\begingroup\makeatletter\ifx\SetFigFont\undefined%
\gdef\SetFigFont#1#2#3#4#5{%
  \reset@font\fontsize{#1}{#2pt}%
  \fontfamily{#3}\fontseries{#4}\fontshape{#5}%
  \selectfont}%
\fi\endgroup%
{\renewcommand{\dashlinestretch}{30}
\begin{picture}(7824,1650)(0,-10)
\put(2712.000,-1950.000){\arc{6000.000}{4.0689}{5.3559}}
\put(5412.000,-1950.000){\arc{6000.000}{4.0689}{5.3559}}
\put(4062.000,-768.750){\arc{3637.500}{3.8759}{5.5488}}
\put(4062.000,-1912.500){\arc{6525.000}{3.9514}{5.4734}}
\path(12,450)(7812,450)
\path(7692.000,420.000)(7812.000,450.000)(7692.000,480.000)
\put(162,0){\makebox(0,0)[lb]{$v_{k}$}}
\put(1062,0){\makebox(0,0)[lb]{$v_{k+1}$}}
\put(1962,0){\makebox(0,0)[lb]{$v_{k+1}$}}
\put(2862,0){\makebox(0,0)[lb]{$v_{k}$}}
\put(3762,0){\makebox(0,0)[lb]{$v_{k}$}}
\put(4662,0){\makebox(0,0)[lb]{$v_{k-1}$}}
\put(5562,0){\makebox(0,0)[lb]{$v_{k}$}}
\put(6462,0){\makebox(0,0)[lb]{$v_{k}$}}
\put(7362,0){\makebox(0,0)[lb]{$v_{k}$}}
\put(1662,1125){\makebox(0,0)[lb]{$B^+_1$}}
\put(3987,1500){\makebox(0,0)[lb]{$I$}}
\put(6462,975){\makebox(0,0)[lb]{$B^-_0$}}
\put(4137,1125){\makebox(0,0)[lb]{$B^-_1$}}
\end{picture}
} }
         	\hspace{-1.9mm}
	\end{array}

$$
we have (assuming $m(s_0)=k$) that
\begin{align*}
\w_\rho(1)&=(1+k)(\l+1-k), & \w_\rho(2)&=\l+2 \\
\w_\rho(3)&=(1+(k-1))(\l+1-(k+1)), & \w_\rho(4)&=-(1+k)(\l+1-k).
\end{align*}
\end{example}



Each coloring $\rho$ of the chords is determined by a subset $V'$ of chords 
such that $\rho(v)=I$ for $v\notin V'$ and by a {\em coloring} $c$ where 
$c(v)$ is an assignement of an element $(\e_v, \d_v) \in \{0,1\}\times 
\{+,-\}$ for each $v\in V'$. Hence we can write 
$$
W_J(D)=J'(\LID(D))|_{x_e=1}
$$
where 
$$
J'(G)=\sum_{V' \subset V} (\l+2)^{|V-V'|}
\sum_{c \, \text{ col of } \, V'} \prod_{v\in V'} \w'_c(v).
$$

An important observation is that $\w'_c(v)$ can be computed in terms of the 
{\em local structure} of the labeled intersection digraph $\LID(D)$. The next 
lemma describes this.

\begin{lemma}
\lbl{lem.2}
Let $\LID(D)=(V,A)$. Then
 $$
 \w'_c(v)= -(-1)^{\e_v} \left(1+\sum_{e \in A, \, v \in e}
z_v(e)\right)\left(\l+1-\sum_{e \in A, \, v \in e} \bar z_v(e)\right)
$$
where $z_v, \bar z_v$ are defined as follows:
\begin{itemize}
\item
If $e$ uncolored then
\item
If $c(v)=(+,\e)$ and $e=(w,v)$ then $z_v(e)=\d_w\e_wx_e$,
\item
If $c(v)=(+,\e)$ and $e=(v,w)$ then $\bar z_v(e)=\d_w\e_wx_e$,
\item
If $c(v)=(-,\e)$ and $e=(w,v)$ then $\bar z_v(e)=\d_w\e_wx_e$,
\item
If $c(v)=(-,\e)$ and $e=(v,w)$ then $z_v(e)=\d_w\e_wx_e$,
\item
If $e$ is red and $e=(v,w)$ then  $z_v(e)=\bar z_v(e)=\d_w\e_wx_e$,
\item
and $z_v(e)=\bar z_v(e)=0$ otherwise.
\end{itemize}
\end{lemma}

\begin{proof}

If $c(v)=(+,\e)$ then $c(v)$ corresponds to the operator $B^+_{\e}$ and hence
$\w'_c(v)|_{x_e=1} = - (-1)^{\e_v} (1+k)(\l+1-k')$. A moment's thought reveals that $k=k_1-k_2+k_3-k_4$ and $k'=k'_1-k'_2+k_3-k_4$ where 
\begin{align*}
k_1 &=|\{e=(w,v) \, \text{uncolored}; \, 
c(w)=(+,1)\}|, & k_2 &=|\{e=(w,v) \, \text{uncolored}; \, c(w)=(-,1)\}| \\
k'_1 &=|\{e=(v,w) \, \text{uncolored}; \,c(w)=(+,1)\}|,  
& k'_2 &=|\{e=(v,w) \, \text{uncolored}; \,c(w)=(-,1)\}|\\
k_3 &=|\{e=(v,w) \, \text{red}; \, c(w)=(+,1)\}|, 
& k_4 &=|\{e=(v,w) \, \text{red}; \,c(w)=(-,1)\}|.
\end{align*}

The reasoning is analogous for $c(v)=(-,\e)$.
\end{proof}

Thus, Theorem \ref{thm.2} follows from the following

\begin{theorem}
\lbl{thm.3}
For all digraphs $G$ with one loop at each vertex, we have
$$J'(G)=J(G).
$$
\end{theorem}

\section{Understanding the state sums $J(G)$ and $J'(G)$}
\lbl{sec.understanding}

In this section we prove Theorem \ref{thm.3}, via a trip to combinatorics
with curious cancellations caused by applications of the 
binomial theorem.

Let us begin by rewriting $J'(G)$. Let $\e_{V'}=\sum_{v \in V'} \e_v$.
Then,
\begin{eqnarray*}
J'(G)&=&
\sum_{V' \subset V} (-1)^{|V'|} (\l+2)^{|V-V'|} \\
& &
\sum_{c \, \text{ col of } \, V'} (-1)^{\e_{V'}}
\sum_{V_1 \subset V', \, V_2 \subset V'}
(-1)^{|V_2|}(\l+1)^{|V'-V_2|}
\prod_{v \in V_1} \left(\sum_{v \in e}
z_v(e)\right)\prod_{v \in V_2}
\left(\sum_{v \in e} \bar z_v(e)\right)
\end{eqnarray*}
where $V_1$ and $V_2$ are possibly overlapping subsets of $V$.

Note that
$$
\prod_{v \in V_1}\left(\sum_{v \in e}z_v(e)\right)
=\sum_f \prod_{v \in V_1}z_v(e(f,v))
$$
where $f:V_1\to A$ maps $v$ to the arc denoted by $e(f,v)$ such that 
$v \in e(f,v)$ and moreover if $e(f,v)$ red then $e=(v,.)$, i.e.
$e$ starts in $v$.
 In other words, $f$ associates with each vertex $v$ of $V_1$ 
an arc incident with it. Similarly, we can rewrite $\prod_{v \in V_2}
\left(\sum_{v \in e}\bar z_v(e)\right)$. Hence,
\begin{eqnarray*}
J'(G)&=&
\sum_{V' \subset V,V_1 \subset V', V_2 \subset V'} 
(-1)^{|V'|} (\l+2)^{|V-V'|}(\l+1)^{|V'-V_2|}(-1)^{|V_2|} \\
& &
\sum_{f:V_1\to A,g:V_2\to A}
\sum_{c \, \text{ col of } \, V'}
(-1)^{\e_{V'}}
\prod_{v \in V_1} \left(
z_v(e(f,v))\right)\prod_{v \in V_2}
\left( \bar z_v(e(g,v))\right).
\end{eqnarray*} 
Let us rewrite the formula more: we fix $W_1=V_1 \cap V_2$, $W_2=V_1 \cup V_2$
and we let $h$ to be disjoint union of $f$ and $g$. 

\begin{remark}
What exactly is $h$? Answer: $h$ is a function that assigns to each vertex of 
$V'$ zero, one or two arcs incident
with it (if the arc is red then it must start in that vertex). Hence 
$|h(v)|\leq 2$ for all $v \in V'$ and if $e \in h(v)$ then $v \in e$.
Here we slowly move towards the formalism of acceptable objects. If $e\in h(v)$
then thicken the end of $e$ containing $v$. Hence $h$ becomes a system of 
thickened
arcs of $G$ so that there are at most two arcs in the system that are thickened
at each vertex of $V'$, and red arcs are thickened always at the start.
\qed
\end{remark}

If we have such an $h$, then $W_1=\{v: |h(v)|=2\}$ and 
$W_2=\{v: |h(v)| \geq 1\}$. Hence, $h$ determines the sets $W_1$ and $W_2$.
Fix an $h$ as above, consider its corresponding sets $W_1,W_2$, and let
$h(W_2)$ denote the system of thickened arcs determined by $h$.  We have
\begin{eqnarray*}
J'(G)&=&
\sum_{V' \subset V} (-1)^{|V'|} (\l+2)^{|V-V'|}
\sum_h A(V',h)
\end{eqnarray*}
where
\begin{eqnarray*}
A(V',h)&=&
\sum_{c \, \text{ col of } \, V'}\sum_{V_2'\subset W_2-W_1}
\sum_{g:W_1 \cup V_2'\to h(W_2): g(v) \in h(v)} B
\end{eqnarray*}
and
\begin{eqnarray*}
B &=&
(-1)^{\e_{V'}}(-1)^{|V_2'\cup W_1|}(\l+1)^{|V'-(W_1 \cup V_2')|}
\prod_{v \in W_1 \cup V_2'} 
\bar z_v(e(g,v))\prod_{v \in W_1 \cup V_1'}
 z_v(e(f,v))
\end{eqnarray*}
where
$V_1'=W_2-(W_1 \cup V_2')$ and $f: W_1 \cup V_1'\to h(W_2)$ is such that
the disjoint union of $f$ and $g$ is $h$.

The next two lemmas restrict the possible configurations of $h$
that contribute non-zero $A(V',h)$.
 
\begin{lemma}
\lbl{lem.ob1}
Let $v \in W_2$. If the only arcs of $h(W_2)$ incident
with $v$ are the arcs of $h(v)$, then $A(V',h)=0$.
\end{lemma}

\begin{proof}
$$
\printname{ob1alt}
	\setlength{\unitlength}{0.02\standardunitlength}
	\begin{array}{c}  \hspace{-1.7mm}
         	\raisebox{-8pt}{\begingroup\makeatletter\ifx\SetFigFont\undefined%
\gdef\SetFigFont#1#2#3#4#5{%
  \reset@font\fontsize{#1}{#2pt}%
  \fontfamily{#3}\fontseries{#4}\fontshape{#5}%
  \selectfont}%
\fi\endgroup%
{\renewcommand{\dashlinestretch}{30}
\begin{picture}(5187,1557)(0,-10)
\path(900,612)(900,12)
\path(3000,12)(3600,612)
\path(4575,612)(5175,12)
\put(3825,1362){\makebox(0,0)[lb]{$v \in W_1$}}
\put(0,1362){\makebox(0,0)[lb]{$v \in W_2-W_1$}}
\path(3600,612)(4050,1062)
\path(3675,537)(3525,687)
\path(4125,987)(3975,1137)
\path(3900,762)(3750,912)
\path(3825,687)(3675,837)
\path(3750,612)(3600,762)
\path(3975,837)(3825,987)
\path(4050,912)(3900,1062)
\path(4125,1062)(4575,612)
\path(4050,987)(4200,1137)
\path(4500,537)(4650,687)
\path(4275,762)(4425,912)
\path(4200,837)(4350,987)
\path(4125,912)(4275,1062)
\path(4350,687)(4500,837)
\path(4425,612)(4575,762)
\path(900,1212)(900,612)
\path(975,1212)(825,1212)
\path(975,612)(825,612)
\path(975,912)(825,912)
\path(975,1062)(825,1062)
\path(975,762)(825,762)
\path(975,1137)(825,1137)
\path(975,987)(825,987)
\path(975,837)(825,837)
\path(975,687)(825,687)
\end{picture}
} }
         	\hspace{-1.9mm}
	\end{array}

$$
Fix $V_2'$, $g$ and $c(V'-\{v\})$. If $B \neq 0$, then the color $\d_v \in \{-,+\}$ of $v$ 
may be determined by $g$ and the orientation of the arcs of $h(v)$. 
However, there is still the choice $\e_v
=0$ or $1$. This influences only $(-1)^{\e_{V'}}$, hence the lemma follows.
\end{proof}

\begin{lemma}
\lbl{lem.ob2}
Let $v \in W_1$. If both arcs of $h(v)$ are uncolored and
oriented in the same way w.r.t. $v$ then $A(V',h)=0$. 
If there are exactly three arcs of $h(W_2)$
incident with $v$ then $A(V',h)=0$. 
\end{lemma}

\begin{proof}
$$
\printname{ob2alt}
	\setlength{\unitlength}{0.02\standardunitlength}
	\begin{array}{c}  \hspace{-1.7mm}
         	\raisebox{-8pt}{\begingroup\makeatletter\ifx\SetFigFont\undefined%
\gdef\SetFigFont#1#2#3#4#5{%
  \reset@font\fontsize{#1}{#2pt}%
  \fontfamily{#3}\fontseries{#4}\fontshape{#5}%
  \selectfont}%
\fi\endgroup%
{\renewcommand{\dashlinestretch}{30}
\begin{picture}(2199,2289)(0,-10)
\path(12,12)(612,612)
\path(1587,612)(2187,12)
\path(1062,2262)(1062,1062)
\put(1287,1887){\makebox(0,0)[lb]{$a$}}
\put(1287,1137){\makebox(0,0)[lb]{$v$}}
\path(612,612)(1062,1062)
\path(687,537)(537,687)
\path(1137,987)(987,1137)
\path(912,762)(762,912)
\path(837,687)(687,837)
\path(762,612)(612,762)
\path(987,837)(837,987)
\path(1062,912)(912,1062)
\path(1137,1062)(1587,612)
\path(1062,987)(1212,1137)
\path(1512,537)(1662,687)
\path(1287,762)(1437,912)
\path(1212,837)(1362,987)
\path(1137,912)(1287,1062)
\path(1362,687)(1512,837)
\path(1437,612)(1587,762)
\end{picture}
} }
         	\hspace{-1.9mm}
	\end{array}

$$
Since $v\in W_1$ we have $|h(v)|=2$. For the first part, if both arcs of $h(v)$ 
are uncolored and oriented in the same way, then there
is no way to choose color $\d_v \in \{-,+\}$ so that $B\neq 0$.
For the second part, if there are exactly three arcs
of $h(W_2)$ incident with $v$ then there is exactly one arc, say $a$, 
which is incident with $v$ and belongs to $h(W_2)-h(v)$.
Fix $V_2'$ and consider pairs $g_1,g_2$ of functions $g$ which differ only on
$v$. If $B \neq 0$ and at least one of the arcs of $h(v)$ is uncolored,
 then the color $\d_v \in \{-,+\}$ of $v$ is determined
by the orientation of the arcs of $h(v)$ and the choice between $g_1$ and 
$g_2$. This color is opposite for $g_1$ and $g_2$,
 and so the edge $a$ is counted with different signs for $g_1$ and $g_2$
while  all the rest remains the same, hence the total contribution is $0$.
If both arcs of $h(v)$ are red then both colors $\d_v \in \{-,+\}$ of $v$
are possible for $B \neq 0$ and both $g_1$ and $g_2$. Hence again
the edge $a$ contributes twice $+1$ and twice $-1$ and the total contribution
is $0$.
\end{proof}

Note that the second property of the above lemma assures that the system
$h(W_2)$ of thickened {\em uncolored} edges is a set for each $h$ which
 contributes a non-zero term to $J'(G)$.

\begin{corollary}
\lbl{cor.ob1}
{\rm (a)} If $e \in h(W_2)$, then both vertices of $e$ belong to $W_2$. \newline
{\em (b)} Each vertex of $W_2$ has degree (i.e., valency) $2$ or $4$ in 
$h(W_2)$ and if $N(v)$ is the set of edges of $h(W_2)$ incident with $v$ then
$|N(v)|=2|h(v)|$.
In other words, the allowed configurations are 
$$
\printname{ob3alt}
	\setlength{\unitlength}{0.02\standardunitlength}
	\begin{array}{c}  \hspace{-1.7mm}
         	\raisebox{-8pt}{\begingroup\makeatletter\ifx\SetFigFont\undefined%
\gdef\SetFigFont#1#2#3#4#5{%
  \reset@font\fontsize{#1}{#2pt}%
  \fontfamily{#3}\fontseries{#4}\fontshape{#5}%
  \selectfont}%
\fi\endgroup%
{\renewcommand{\dashlinestretch}{30}
\begin{picture}(4299,2214)(0,-10)
\path(87,912)(87,12)
\path(2112,12)(2712,612)
\path(3687,612)(4287,12)
\path(2712,612)(3162,1062)
\path(2787,537)(2637,687)
\path(3237,987)(3087,1137)
\path(3012,762)(2862,912)
\path(2937,687)(2787,837)
\path(2862,612)(2712,762)
\path(3087,837)(2937,987)
\path(3162,912)(3012,1062)
\path(3237,1062)(3687,612)
\path(3162,987)(3312,1137)
\path(3612,537)(3762,687)
\path(3387,762)(3537,912)
\path(3312,837)(3462,987)
\path(3237,912)(3387,1062)
\path(3462,687)(3612,837)
\path(3537,612)(3687,762)
\path(87,2112)(87,1512)
\path(162,2112)(12,2112)
\path(162,1512)(12,1512)
\path(162,1812)(12,1812)
\path(162,1962)(12,1962)
\path(162,1662)(12,1662)
\path(162,2037)(12,2037)
\path(162,1887)(12,1887)
\path(162,1737)(12,1737)
\path(162,1587)(12,1587)
\path(87,1512)(87,912)
\path(162,1512)(12,1512)
\path(162,912)(12,912)
\path(162,1212)(12,1212)
\path(162,1362)(12,1362)
\path(162,1062)(12,1062)
\path(162,1437)(12,1437)
\path(162,1287)(12,1287)
\path(162,1137)(12,1137)
\path(162,987)(12,987)
\path(3162,1062)(4287,2187)
\path(3237,1062)(2112,2187)
\end{picture}
} }
         	\hspace{-1.9mm}
	\end{array}

$$
\end{corollary}

\begin{proof}
It follows from Lemmas \ref{lem.ob1} and \ref{lem.ob2} that $|N(v)| \geq
2|h(v)|$ at each vertex $v$. On the other hand, each edge of $h(W_2)$ has one
{\em thick end} and one {\em thin end}, and so there cannot be more thin
ends than thick ends. Hence $|N(v)| =2|h(v)|$ and the corollary follows.
\end{proof}

\begin{corollary}
\lbl{cor.ob2}
If $W_2 \neq V'$, then $A(V',h)=0$.
\end{corollary}

\begin{proof}
We write
$$
A(V',h)=\sum_{c \, \text{ col of } \, V'-W_2}
(-1)^{\e_{V'-W_2}} \,\, (\text{rest})
$$
where the 'rest' is not influenced by the colorings in $V'-W_2$. Hence,
$$
A(V',h)= (\text{rest}) \, \sum_{C \subset V'-W_2} (-1)^{|C|} 
2^{|V'-W_2|}
$$
which vanishes unless $V'=W_2$.
\end{proof}

Summarizing, a function $h$ such that $A(V',h)\neq 0$ determines a collection
of thickened arcs that is almost an acceptable object:
\begin{itemize}
\item
 each vertex has degree $2$ or $4$ in $V'$ and $0$ 
in $V-V'$,
\item
exactly half of the arcs incident with a vertex are thickened
at that vertex,
\item
if there are two uncolored arcs thickened at a vertex then they
have opposite orientation with respect to the vertex,
\item
the red arcs are always thickened at the start.
\end{itemize}
 Let us call such  object {\em good} on $V'$.
Note also that for each coloring of a good object on $V'$ 
which contributes non-zero to $B$ we must have $\e_v=1$ for each $v\in V'$ and hence
$(-1)^{\e_{V'}}=(-1)^{|V'|}$. Collecting
our efforts so far, we have
\begin{equation}
\lbl{eq.1}
J'(G) = \sum_{V' \subset V}  
(\l+2)^{|V-V'|} \sum_{K \,\text{ good on }\, V'} A'(V',K)
\end{equation}
where
$$
A'(V',K)=\sum_{V_2' \subset V'-W_1}
(\l+1)^{|V'-(W_1 \cup V_2')|} C(V_2',V',K),
$$
$W_1$ is the set of vertices of $V'$ of degree $4$ in $K$ and 
$$
C(V_2',V',K)=\sum_{g':W_1\to K} (-1)^{|W_1|}(-1)^{|V_2'|}
\sum_{c \, \text{ col }}
\prod_{v \in W_1 \cup V_2'} 
\bar z_v(e(g,v))\prod_{v \in V'-V_2'}
 z_v(e(f,v)),
$$
and $g', g, f$ have the following properties: 
\begin{itemize}
\item if $p\in \{g',g,f\}$ then $p(v)$ is an arc of $K$ incident 
with $v$ and thickened at $v$,
\item $g:W_1\cup V_2'\to K$ is unique such function extending $g'$,
\item $f: V'-V_2' \to K$ is unique such function with $f\cup g=K$.
\end{itemize}

\begin{lemma}
\lbl{lem.ob3}
$C(X,V',K)=C(Y,V',K)$ for arbitrary $X, Y$ subsets of $V'-W_1$.
\end{lemma}

\begin{proof}
Exactly half of the edges of $K$ incident with each vertex are thickened 
at that vertex and hence $K$ may be regarded as a 
union of cycles $Z_1,\dots,Z_m$ such that each $Z_i$ has the form
$$
\psdraw{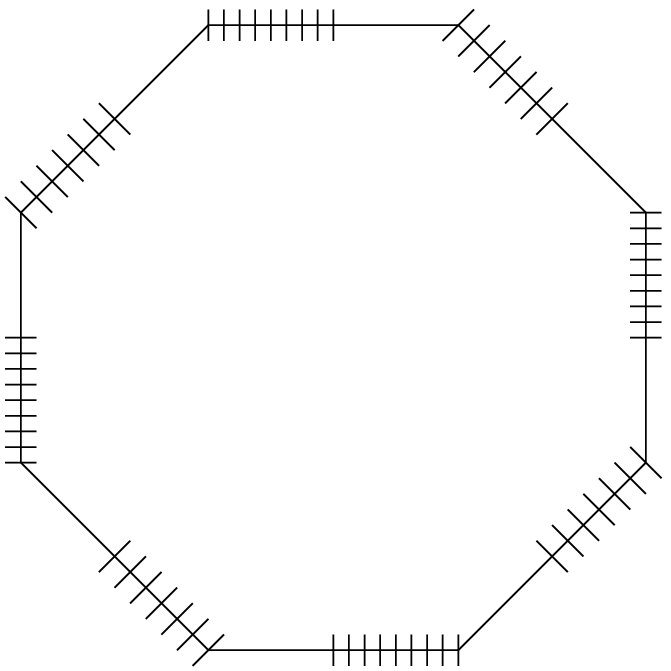}{1in}
$$
and such that each vertex of $K$ lies in at most two of the cycles $Z_i$.
In other words, we may think that $K$ is pictured schematically as follows
$$
\psdraw{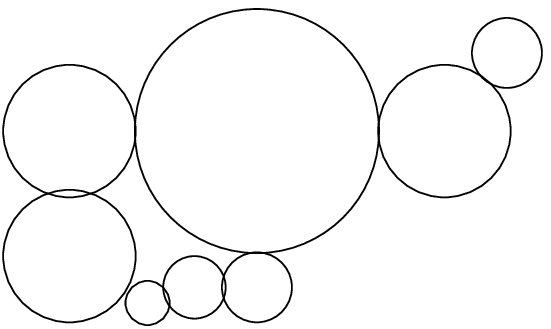}{1in}
$$
(where for simplicity, we have drawn the cycles $Z_i$ as circles). 
As we observed above, 
$$
C(X,V',K)=\sum_{g':W_1\to K} (-1)^{|W_1|}(-1)^{|X|}
\sum_{c \, \text{ col }}
\prod_{v \in W_1 \cup X} 
\bar z_v(e(g,v))\prod_{v \in V'-X}
 z_v(e(f,v))
$$
It suffices to show the following \newline
{\bf Claim.}
$$
C(X,V',K)=\sum_{K'}(-1)^{a(K')}x_{K}
$$ 
where $K'$ is any collection of thickened arcs obtained from $K$
by changing  orientation of some (possibly none) red arcs of $K$
 so that $K'$ is an acceptable object (i.e. if two arcs are thickened
at a vertex then they are oppositely oriented w.r.t. that vertex),
and $a(K')$ is the number of arcs of $K'$ directed in agreement 
with the thickening.

\noindent
{\em Proof of the Claim.}
We can write 
$$
C(X,V',K)=\sum_{g':W_1\to K} D
$$
and
$$
D=(-1)^{|W_1|}(-1)^{|X|}
\sum_{c \, \text{ col }}
\prod_{v \in W_1 \cup X} 
\bar z_v(e(g,v))\prod_{v \in V'-X}
 z_v(e(f,v)).
$$
First we observe that the Claim is true when $K$ has a vertex
of degree 2 where the thickened edge is red. Indeed, in this case
both sides of the formula in the Claim equal $0$.
Hence let $K$ donot have such a vertex.
Let $S$ be the set of vertices of $K$ where two red arcs are thickened.
Observe that exactly $2^{|S|}$ colorings $c$ contribute a non-zero term
to $D$: 
we observed before that necessarily $\e(v)=1$ for each $v\in V'$.
Moreover, each color $\d(v)$ contributing non-zero to $D$ is uniquely
determined for each $v$ where at least one uncolored arc is thickened:
 explicitly, let $v$ be a vertex
of $Z_i$ and let $e$ be the unique arc of $Z_i$ thickened at $v$, and
let $e$ be uncolored.
Then $\d(v)$ depends on the orientation of $e$ and whether $e$ belongs to $g$
or not. Hence there is at most one coloring $\d(v)$
which contributes a non-zero term to $D$. Also observe that the unique
'non-zero coloring' of these vertices of each $Z_i$ compose well together:
this follows from the fact that if $v$ is a vertex of degree $4$ in $K$
then exactly $2$ arcs are thickened at $v$, and if both are uncolored
then they have different orientation
with respect to $v$ and one belongs to $g$ and the other belongs to $f$.
Finally observe that for vertex $v$ where two red arcs are thickened,
any $\d(v)$ contributes a non-zero.

Next we will  observe that the contribution of each of these
$2^{|S|}$ colorings to $D$ is the same and equals $(-1)^{a(K')}x_{K}$ 
where $K'$ is any collection of thickened arcs obtained from $K$
by changing  orientation of some (possibly none) red arcs of $K$
so that $K'$ is an acceptable object (i.e., if two arcs are thickened
at a vertex then they are opositely oriented with respect to that vertex).
This proves the Claim since the number of objects $K'$ equals
$2^{|S|}$. 
Hence it remains to confirm the contribution of each of the allowed
colorings to $D$.

First observe that it is true when $X =\emptyset$.
Next, let us put a vertex $v$ of $Z_i-W_1$ into $X$ and let $e$ be the arc
of $Z_i$ thickened at $v$. Then $e$ is uncolored by our assumption and
we need to change $\d(v)$ in order to have a nonzero contribution, hence the product
of signs along $Z_i$ changes but $(-1)^{|X|}$
also changes and so the final total sign is the same.


\end{proof}

We let 
\begin{eqnarray*}
C(V',K) &=& \sum_{g':W_1\to K} \sum_{K'}(-1)^{a(K')}x_{K} \\
&=&
2^{\deg_4(K)}\sum_{K'}(-1)^{a(K')}x_{K}.
\end{eqnarray*}

Equation \eqref{eq.1} together with Lemma \ref{lem.ob3} implies that
\begin{eqnarray*}
J'(G) &=& 
\sum_{V' \subset V}  
(\l+2)^{|V-V'|} \sum_{K \,\text{ good on }\, V'} C(V',K)
\sum_{A \subset V'-W_1} (\l+1)^{|A|} \\
& = &
\sum_{V' \subset V}  
(\l+2)^{|V-V'|} \sum_{K \,\text{ good on }\, V'} C(V',K)
(\l+2)^{|V'-W_1|} \\
& = &
\sum_{K \, \text{good on}\, V'} 
(\l+2)^{|V-W_1|} C(V',K) \\
& = &
\sum_{K \, \text{acceptable}} (\l+2)^{|V|-\deg_4(K)} 2^{\deg_4(K)} (-1)^{a(K)}x_K
\\
& = & J(G)
\end{eqnarray*}
which concludes the proof of Theorem \ref{thm.3}.
\qed

\section{Converting to Permanents}
\lbl{sec.per}

The goal of this section is to convert the state sum $J(D)$ into
a permanent, in the following way:

\begin{theorem}
\lbl{thm.4}
$$
J(\LID(D))|_{x_e=1}=\Per(\IMJ(D))
$$
\end{theorem}

Note that Theorem \ref{thm.1} follows from Theorems \ref{thm.2} and 
\ref{thm.4}.

We will achieve the conversion of $J(\LID(D))$ in a permanent by
a local modification (we can say, a blow-up) of each of the vertices of the 
original digraph $\LID(D)$.
It turns out that the modification triples each of the vertices of $G$.
Why triple? Because in a sense $J(\LID(D))$ has to do with the $\fsl_2$-Lie algebra.
Thus, we are back to Lie algebras, this time through a common blow-up
trick of combinatorics.
 
\noindent
{\em Proof (of Theorem \ref{thm.4}).}
We have that
$$
J(\LID(D))|_{x_e=1}=\sum_{K \, \text{acceptable}} (\l+2)^{|V|-\deg_4(K)} 
2^{\deg_4(K)} (-1)^{a(K)},
$$
where $a(K)$ equals the number of arcs of $K$ thickened in agreement with 
their orientation. A single loop of $K$ is always leaving its vertex, and it is
directed in agreement with its thickening. Hence each single loop contributes
'(-1)' to $a(K)$. Let us now get rid of these single loops: an acceptable 
object
without single loops will be called {\em connected} and a connected object
where each vertex of $V$ has degree at least $2$ will be {\em super}. We have
\begin{eqnarray*}
J(\LID(D))|_{x_e=1}&=& \sum_{K \, \text{acceptable}} (\l+2)^{|V|-\deg_4(K)} 
2^{\deg_4(K)} (-1)^{a(K)} \\
&=&
\sum_{K \, \text{connected}} (\l+2)^{|V|-\deg_4(K)} 
2^{\deg_4(K)} (-1)^{a(K)}\sum_{U\subset V-K} (-1)^{|U|} \\
&=& 
\sum_{K \, \text{super}} (\l+2)^{|V|-\deg_4(K)} 2^{\deg_4(K)} (-1)^{a(K)}.
\end{eqnarray*}
 
We will use a variation $\TID(D))$, the {\em thickened intersection digraph} 
of $D$ defined as follows. Double each of the arcs of $\LID(D)$ and:
\begin{enumerate}
\item
 if the arc is uncolored then thicken such pair at opposite ends,
\item
if the arc is red then thicken each arc of the pair at the start,
and then change the orientation of one of them,
\item
 thicken each loop at its initial segment, i.e., in agreement with 
its orientation. 
\end{enumerate}
In pictures, the thickening of $\LID$ is the substitution
$$
\psdraw{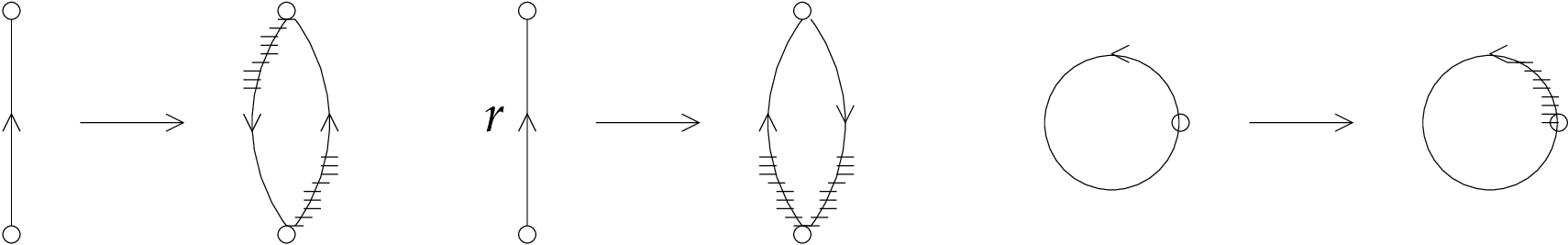}{4in}
$$


Now we can write 
$$
J(\LID(D))|_{x_e=1}= \sum_{K \, \text{super subobject of }\TID(D)}  
(\l+2)^{|V|-\deg_4(K)} 2^{\deg_4(K)} (-1)^{a(K)}.
$$
 
 
Let us describe now how a thickened digraph $D(\IMJ)$ may be constructed
 from $\TID(D)$. The 
construction easily follows from the definition of matrix $\IMJ$: it
consists in replacing each vertex of $\TID(D)$ by a 'gadget' on three vertices,
as follows: 
$$ 
\psdraw{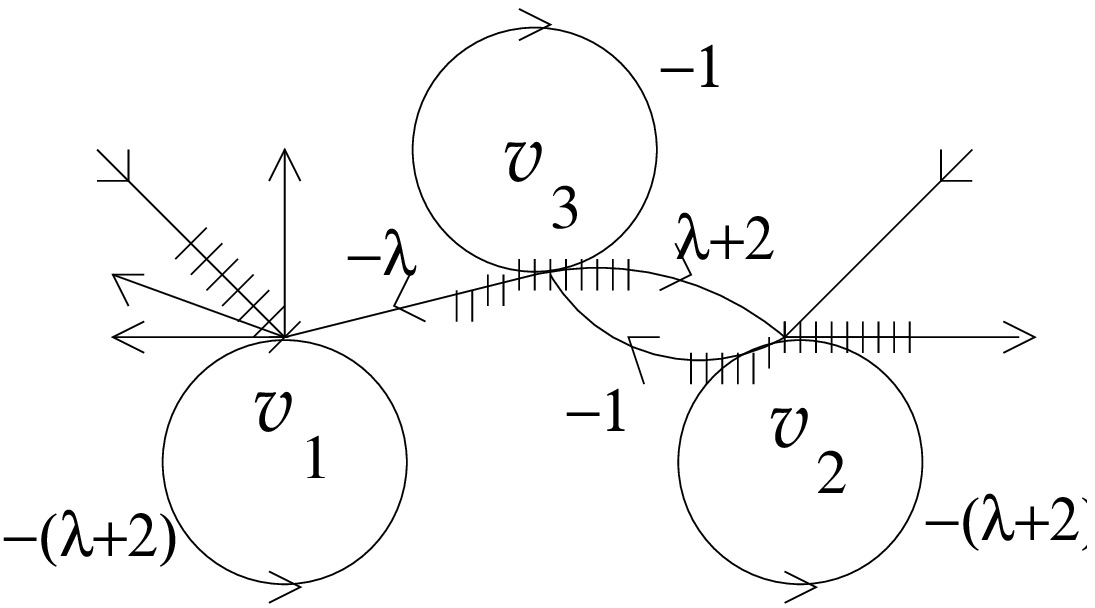}{2in} 
$$
The construction goes as follows:
\begin{enumerate}
\item[1.] For each vertex $v$ of $\TID(D)$ introduce three vertices 
$v_1,v_2,v_3$ for $D(\IMJ)$.
\item[2.] Define the thickened arcs and their weights among each triple 
$v_1,v_2,v_3$ as follows:
\begin{itemize}
\item add loop $l_i$ at each $v_i$ and let $w(l_1)=w(l_2)=-(\l+2)$ and 
$w(l_3)=-1$,
\item add arc $(v_3,v_1)$ thickened at $v_3$ with weight $-\l$,
\item add arc $(v_3,v_2)$ thickened at $v_3$ with weight $\l+2$,
\item add arc $(v_2,v_3)$ thickened at $v_2$ with weight $-1$.
\end{itemize}
\item[3.] For each thickened arc $(u,w)$ of $\TID(D)$ do the following:
\begin{itemize}
\item If $(u,w)$ is thickened at $u$ then add $(u_2,w_1),(u_2,w_2)$ thickened 
at $u_2$ with weights equal to $1$,
\item If $(u,w)$ is thickened at $w$ then add $(u_1,w_1),(u_2,w_1)$ thickened 
at $w_1$ with weights equal to $1$.
\end{itemize}
\end{enumerate}

It follows directly from the definition of the permanent that
$$
\Per(\IMJ)=\sum_{L\in \cL} (-1)^{a(L)} w_L,
$$
where $\cL$ is the set of all acceptable subobjects of $D(\IMJ)$ where each 
degree equals $2$.

We need to show that
\begin{eqnarray*}
\sum_{\begin{subarray}{c} L \, \text{with each degree 2} \\
     \text{ acceptable subobject of} D(\IMJ)
      \end{subarray}} (-1)^{a(L)} w_L &=&
\sum_{K \, \text{super subobject of }\TID(D)}  (\l+2)^{|V|-\deg_4(K)} 
2^{\deg_4(K)} (-1)^{a(K)}.
\end{eqnarray*}
We will prove it by constructing a partition of acceptable subobjects of 
$D(\IMJ)$ where each vertex has degree 2, and associating each partition 
class which contributes non-zero to $\Per(D(\IMJ))$ with uniquelly determined 
super subobject of $\TID(D)$.
 
Let $L$ be an acceptable subobject of $D(\IMJ)$ where each vertex
has degree 2. Denote by $OL$ the set of all thickened arcs of type $(u_iw_j)$,
$u\neq w$, and let $IL=L-OL$. Note that if we forget the lower indices
at vertices, $OL$ naturally 
corresponds to a set $OK$ of thickened arcs of $\TID(D)$. Note that no arc
in $OK$ is a loop and $a(OL)=a(OK)$. It remains to be seen what to do with the
thickened arcs of $IL$. We may consider each triple of vertices $v_1,v_2,v_3$ 
separately. Let $ILv$ denote the set of thickened arcs of $L$ among 
$v_1,v_2,v_3$. We distinguish four cases.
\begin{itemize}
\item  
$ILv$ consists of all three loops or the loop at $v_1$ and the arcs 
$(v_2,v_3),(v_3,v_2)$.
Let $C_0$ be the class of all $L$ which in at least one triple $v_1,v_2,v_3$ 
behave in this way. Note that the total contribution of $C_0$ to $\Per(D(\IMJ))$ 
is 0, and so we may assume that this case never happens (it corresponds to 
single loop at $v$ in $\TID(D)$ which is not allowed for supersubobjects).
\item
$ILv$ consists of loop at $v_1$ or loop at $v_2$ (but not both), and loop at 
$v_3$. Then let $IKv=\emptyset$. Hence $v$ will have degree 2 in $K$ and 
contribute $(\l+2)$ to $J(\LID(D))$, which exactly equals to the contribution 
of $ILv$ to $\Per(D(\IMJ))$.
\item
$ILv$ consists of $\{(v_2,v_3), (v_3,v_2)\}$ or  $\{(v_3,v_1), (v_2,v_3)\}$. 
Then $v$ has degree $2$ in $OK$ and the edge of $OK$ incident and thickened at 
$v$ is entering $v$. In this case let $IKv$ consist of the loop at $v$. Note 
that the total contribution of $ILv$ to $\Per(D(\IMJ))$ is
$-(\l+2)+\l=-2$ and $2(-1)$ is also the contribution of $IKv$ to $J(\LID(D))$.
\item
$ILv=\emptyset$. Then $v$ has degree 4 in $OK$ and we let $IKv=\emptyset$. In 
this case each $v_i$, $i=1,2$ is incident with one arc of $L$ not thickened 
at $v_i$. Let $L'$ be obtained from $L$ by exchanging the incidence of 
non-thickened arcs between $v_1$ and $v_2$. The total contribution of $ILv$ 
and $IL'v$ to $\Per(D(\IMJ))$ is
$2$ and $2$ is also the contribution of $IKv$ to $J(\LID(D))$.
\end{itemize}
It is easy to check that we indeed partitioned the set of all acceptable 
subobjects of $D(\IMJ)$ where each vertex has degree 2 and that we exhausted 
all super subobjects of $\TID(D)$. This finishes the proof of Theorem 
\ref{thm.4}.

\section{Proof of the corollaries}
\lbl{sec.cor}

Corollary \ref{cor.2} is immediate. For Corollary \ref{cor.3}, observe
that $W_{JJ}=W_{JJ^{(0)}}$ is given by Corollary \ref{cor.2} for $n=0$. We 
claim that this formula equals to $\Per(\IM)$. This may be observed as 
follows: assume $i < j$. Associate $\IM(D)_{ij}$ with the
arc $(i,j)$ of $\LID(D)$ thickened at $i$ and $\IM(D)_{ji}$ with arc $(i,j)$ 
thickened at $j$. This associates, with each term of the expansion of 
$\Per(\IM)$, an acceptable object of uncolored arcs only,
 with each degree equal to 2 and no loops. 
Denote the set of such acceptable objects by $\K_1$. It is straightforward to 
check that 
$$
\Per(\IM)=\sum_{K\in \K_1} (-1)^{a(K)}.
$$
On the other hand, $W_{JJ^{(0)}} = \sum_{K\in \K_2} (-1)^{a(K)}$, where $\K_2$
is the set of all acceptable objects where each degree is 0 or 2 (loops are 
allowed and they contribute 2 to the degree).
First observe that the contribution of the acceptable objects of $\K_2$ 
that contain a red arc cancels out since we can change the orientation
of a red arc in such an object, and get again an object of $\K_2$,
with oposite contribution. Hence assume $\K_2$ has no objects with
red arcs. If $K \in \K_2$ then let $L(K)$ 
denote the acceptable subobject of $K$ obtained from $K$ by deleting all loops.
If $L=L(K)$ let $V(L)$ denote the set of vertices of $L$ of non-zero degree 
and let $\E(L)=\{K\in \K_2; L(K)=L\}$. By the binomial theorem,
$$
\sum_{K\in \E (L)} (-1)^{a(K)}= (-1)^{a(L)} \sum_{W \subset (V-V(L))} 
(-1)^{|W|}=0
$$
whenever $V(L) \neq V$. This proves the corollary.

\section{Numerical examples}
\lbl{sec.numerical}

\subsection{Running the program}

The formula of Theorem \ref{thm.1} for the weight system of the colored Jones
function is easy to program, as is demonstrated by the source code given in 
the next section. To run the program, first start {\tt Mathematica} \cite{Wo}
and load the \verb$JonesPermanent$ package of the next section (available 
from our web-site, too):

{\small
\begin{verbatim}
Mathematica 4.1 for Linux
Copyright 1988-2000 Wolfram Research, Inc.
 -- Motif graphics initialized --
In[1]:=<< JonesPermanent.m
\end{verbatim}
}

\vskip 2mm
For the chord diagram \verb$CDP[5,7,6,8,1,3,2,4]$ of Example \ref{ex.1},
$W_J$ is given by
\vskip 2mm

{\small
\begin{verbatim}
In[2]:=formula[CDP[5,7,6,8,1,3,2,4]]
Out[2]:=4 x^4 + 16 x^3 -4 x^2 -40 x
\end{verbatim}
}

\vskip 2mm
where $\l=x$, and the notation \verb$CDP[a_1,...,a_{2n}]$ for a permutation
$\s=(a_1,\dots,a_{2n})$ of $(1,\dots,2n)$ of order two means the chord diagram 
of the $2$-cycles of $\s$.


\subsection{The source code}

{\small 
\begin{verbatim}
CD2CDP[cd_CD] := 
  Module[{m = Length[cd]/2, i, beg, end}, 
   beg = Table[Position[cd, i], {i, 1, m}]; end = Reverse /@ beg; 
    Apply[CDP, Range[2*m] /. Thread[Flatten[beg] -> Flatten[end]]]]

An:={{0,0,0},{-1,-1,0},{0,0,0}}; Ap:={{1,1,0},{0,0,0},{0,0,0}};
Ar:={{1,1,0},{-1,-1,0},{0,0,0}}; Az:={{x+2,0,0},{0,x+2,1},{x,-x-2,1}}

TheSymbol[a_,b_]:=
Module[{a1,a2,b1,b2},   (* a,b are chords *)
      a1=a[[1]]; a2=a[[2]]; b1=b[[1]]; b2=b[[2]];
      If[a==b, Zero, 
      If[a1 < b1 < a2 < b2, Minus, 
      If[b1 < a1 < b2 < a2, Plus,
      If[b1 < a1 < a2 < b2, Red, Nought]]]]]

rule[a_,b_,P_] :=
Module[{a1,a2,b1,b2,ena,duo},  (* a,b are numbers *)
	 a1=Ceiling[a/3]; a2=If[a==3*Floor[a/3], 3, a-3*Floor[a/3]];
         b1=Ceiling[b/3]; b2=If[b==3*Floor[b/3], 3, b-3*Floor[b/3]];
         ena=P[[a1]]; duo=P[[b1]];
         If[TheSymbol[ena,duo]===Zero, Az[[a2,b2]],
           If[TheSymbol[ena,duo]===Minus, An[[a2,b2]],
             If[TheSymbol[ena,duo]===Plus, Ap[[a2,b2]],
               If[TheSymbol[ena,duo]===Red, Ar[[a2,b2]], 0]]]]]

Permanent[m_List] := With[{v = Array[x, Length[m]]},
                           Coefficient[Times @@ (m.v), Times @@ v] ] 

formula[cdp_CDP] := 
   Module[{m, chords, k, kp, km, e, g, i, n2, j1, j2}, 
     m = Length[cdp]/2; chords = 
      Apply[CL, Union[Sort /@ Transpose[{Range[2*m], Apply[List, cdp]}]]];
     Expand[Permanent[Table[rule[j1,j2,chords],{j1,3*m},{j2,3*m}]]]] 

 
\end{verbatim}}

\ifx\undefined\bysame
	\newcommand{\bysame}{\leavevmode\hbox
to3em{\hrulefill}\,}
\fi

\end{document}